# RESEARCH ANNOUNCEMENT

# WIENER'S TAUBERIAN THEOREM IN $L^1(G//K)$ AND HARMONIC FUNCTIONS IN THE UNIT DISK

Y. BEN NATAN, Y. BENYAMINI, H. HEDENMALM, AND Y. WEIT

ABSTRACT. Our main result is to give necessary and sufficient conditions, in terms of Fourier transforms, on a closed ideal $I$ in $L^1(G//K)$, the space of radial integrable functions on $G = SU(1,1)$, so that $I = L^1(G//K)$ or $I = L^1_0(G//K)$—the ideal of $L^1(G//K)$ functions whose integral is zero. This is then used to prove a generalization of Furstenberg's theorem which characterizes harmonic functions on the unit disk by a mean value property and a "two circles" Morera type theorem (earlier announced by Agranovskiĭ).

## 1. INTRODUCTION

Let $G$ be a locally compact abelian group. By Wiener's Tauberian theorem, if the Fourier transforms of the elements of a closed ideal $I$ of $L^1(G)$ do not have a common zero, then $I = L^1(G)$.

In the non-abelian case the analogue of Wiener's theorem for two-sided ideals holds for all connected nilpotent Lie groups and all semi-direct products of abelian groups [16]. However, Wiener's theorem does not hold for any non-compact connected semisimple Lie group [7], [16].

In their seminal series of papers on harmonic analysis on $SU(1,1)$, Ehrenpreis and Mautner use the ideal structure in the disk algebra $A(\mathbb{D})$ to show that the analogue of Wiener's theorem fails even for the commutative subalgebra $L^1(G//K)$ of spherical functions ([7], see also [3]). They realized that in addition to the non-vanishing of the Fourier transforms, a condition on the rate of decay of the Fourier transforms at infinity is also necessary. For technical reasons, it was necessary

Received by the editors May 31, 1994.
1991 *Mathematics Subject Classification.* Primary 43A90.
*Key words and phrases.* Spherical functions, $SU(1,1)$, Resolvent transform, Spectral synthesis, Wiener's theorem, $\mu$-harmonic functions, two circles theorems.

This is a part of the first author's Ph.D. dissertation, prepared at the Technion - Israel Institute of Technology under the supervision of the second and the fourth authors.

The second author's work was partially supported by the fund for the promotion of research at the Technion - Israel Institute of Technology.

The third author's work was partially supported by the Swedish Natural Science Research Council and by the Wallenberg Prize from the Swedish Mathematical Society.

The final detailed version of this paper will be submitted for publication elsewhere.

for them to impose various smoothness conditions (of the Fourier transforms) in addition to the natural conditions of non-vanishing of the Fourier transforms and the "correct" rate of decay, in their analogue of Wiener's theorem ([6], see also [3]).

It is known that smoothness conditions make Wiener's theorem much easier. See, for example, [12] for a trivial proof that if $f \in L^1(\mathbb{R})$ and its Fourier transform $\widehat{f}$ is twice differentiable and never vanishes, then the closure of the ideal generated by $f$ is all of $L^1(\mathbb{R})$.

The main result of the present paper is *a genuine analogue of Wiener's theorem without any superfluous smoothness condition.*

We use the method of the resolvent transform, as developed by Gelfand, Beurling, and Carleman. Gelfand's point of view was later rediscovered and elaborated upon by Domar and was used by Hedenmalm and others in the study of closed ideals in Banach algebras [4], [5], [10], [2].

As applications of the "correct" version of Wiener's theorem, we follow the ideas of [3] and give a generalization of a theorem of Furstenberg [8, 9] characterizing bounded harmonic functions in the unit disk as the bounded solutions of certain convolution equations (i.e., $\mu$-harmonic functions) and a "two circles" Morera-type theorem characterizing holomorphic functions in the unit disk.

## 2. The main result

Let $G = SL_2(\mathbb{R})$ be the group of all $2 \times 2$ real matrices with determinant 1. $G$ is equivalent to $SU(1,1)$ and to the group of all conformal automorphisms of the unit disk $\mathbb{D}$. $G$ has a polar decomposition $G = KA^+K$, where $K$ is the subgroup of all rotations of $\mathbb{D}$ and $A^+$ is identified with $\mathbb{R}^+$. The Haar measure on $G$ is $dg = \sinh 2\zeta \, d\zeta \, d\varphi \, d\theta$, where $d\zeta$ is the Lebesgue measure on $\mathbb{R}^+$, and $d\varphi$ and $d\theta$ are the Haar measure on the unit circle $K$. The symmetric space $G/K$ is identified with *the Poincaré model of the hyperbolic plane* $\mathbf{H}^2$, i.e., with $\mathbb{D}$. The space $K\backslash(G/K)$ which is denoted by $G//K$ is a measure space with the measure $\sinh 2\zeta \, d\zeta$. Changing variables by $x = \cosh 2\zeta$, we can use the Lebesgue measure $\frac{1}{2}dx$ on $\mathbb{R}^+$. The algebra $L^1(G//K)$ is a commutative Banach algebra, with the maximal ideals space identified with the quotient of $\mathcal{S} = \{s \in \mathbb{C} : 0 \leq \Re s \leq 1\}$ under the identification $s = 1 - s$. The Gelfand transform $\widehat{f}(s)$ of $f \in L^1(G//K)$ is given by $\widehat{f}(s) = \frac{1}{2}\int_1^\infty f(x)P_{s-1}(x)\,dx$ where $s \in \mathcal{S}$, and $P_{s-1}$ is the Legendre function of the first kind. It is continuous in the strip $\mathcal{S}$ and is analytic in its interior.

We refer to [11] and [14] for details and to [15], [17], and [18] for properties of special functions.

In analogy to the Beurling-Rudin description of primary closed ideals in $A(\mathbb{D})$, we study rates of decay of Fourier transforms of elements of ideals in $L^1(G//K)$. Following [10, p. 133], we define for $f \in L^1(G//K)$

$$\delta_\infty(f) = -\limsup_{t \to +\infty} e^{-\pi t} \log |\widehat{f}(\tfrac{1}{2} + it)|,$$

$$\delta_0(f) = -\limsup_{x > 0,\, x \to 0} x \log |\widehat{f}(x)|.$$

The main result is the analogue in $L^1(G//K)$ of Wiener's Tauberian theorem with exact decay condition and without superfluous technicalities. Recall that the hull of an ideal is the set of common zeros of the Fourier transforms of its elements.

**Theorem 1.** *Let $M \subseteq L^1(G//K)$ generate the closed ideal $I(M)$ in $L^1(G//K)$. Assume $\inf_{f \in M} \delta_\infty(f) = 0$.*

*(a) If $I(M)$ has an empty hull, then $I(M) = L^1(G//K)$.*

*(b) If $\operatorname{hull}(I(M)) = \{0, 1\}$ and $\inf_{f \in M} \delta_0(f) = 0$, then $I(M) = L_0^1(G//K)$, the space of functions in $L^1(G//K)$ whose integral vanishes.*

## 3. Applications

Let $f$ be harmonic on $\mathbb{D}$. For any $g \in G$, $f \circ g$ is also harmonic, and its value at zero is the average of its values on any circle centered at $0$. Identifying points in $\mathbb{D}$ with elements of $G/K$ and radial measures on $\mathbb{D}$ with bi-invariant measures on $G$, we obtain that if $\mu$ is a radial measure on $\mathbb{D}$ with $\mu(\mathbb{D}) = 1$, then

$$\int_G f(gh)\, d\mu(h) = f(g) \qquad \text{for all } g \in G$$

whenever $f \in L^\infty(G//K)$ is harmonic on $\mathbb{D}$, i.e., $f$ is $\mu$-harmonic [8]. Conversely, one may ask under what conditions on a radial measure $\mu$ are the only $\mu$-harmonic functions $f \in L^\infty(G//K)$ harmonic? Using Theorem 1, we obtain

**Theorem 2.** *Let $\mu$ be a radial measure on $\mathbb{D}$ such that $\mu(\mathbb{D}) = 1$, $\mu$ has no atom at $0$, $\hat{\mu}(s) \neq 1$ for $s \in \mathcal{S}$, $s \neq 0, 1$, and so that*

$$\limsup_{x>0,\, x\to 0} x \log|1 - \hat{\mu}(x)| = 0.$$

*Then every bounded $\mu$-harmonic function is harmonic.*

Theorem 2 generalizes the following result of Furstenberg [9].

**Corollary 3.** *If $\mu$ is a radial probability measure on $\mathbb{D}$ with no atom at 0, then every bounded $\mu$-harmonic function is harmonic.*

We finish with a "two circles" Morera-type characterization of holomorphic functions in the unit disk announced earlier by Agranovskiĭ.

Fix two central circles $\gamma_1$ and $\gamma_2$ of radii $r_1$ and $r_2$ respectively, and define for $j = 1, 2$

$$J_j(s) = {}_2F_1(2 - si, \frac{3}{2}; 3; -4r_j(1 - r_j)^{-2})$$

where ${}_2F_1$ is the hypergeometric function.

**Theorem 4.** *If $f$ is a measurable function on $\mathbb{D}$ satisfying the growth condition $|f(z)| \leq c(1 - |z|^2)^{-1}$ in $\mathbb{D}$, so that for $j = 1$ and $2$*

$$\int_{g(\gamma_j)} f(z)\, dz = 0, \quad \text{for almost all } g \in G,$$

*and if $J_1(s)$ and $J_2(s)$ have no common zero in the strip $\mathcal{S}$, then $f$ coincides almost everywhere in $\mathbb{D}$ with a holomorphic function.*

## 4. Sketches of proofs

*Proof of Theorem 1.* The proof of this theorem is carried out in eight steps.

(1) For every $\lambda \in \mathbb{C} \setminus \mathcal{S}$ find $b_\lambda \in L^1(G//K)$ such that

$$\hat{b}_\lambda(s) = (s(1 - s) - \lambda(1 - \lambda))^{-1}.$$

(2) For every $g \in L^\infty(G//K)$ which annihilates $I(M)$ define the resolvent transform $\mathfrak{C}[g]$ by $\mathfrak{C}[g](\lambda) = \langle b_\lambda, g \rangle$ for $\lambda \in \mathbb{C} \setminus \mathcal{S}$.

(3) Show that $\mathfrak{C}[g]$ is a holomorphic function outside $\mathcal{S}$ which may be extended to an entire function in case (a) and to a holomorphic function in $\mathbb{C} \setminus \{0, 1\}$ in case (b).

(4) Show that $\mathfrak{C}[g](\lambda) = \langle \frac{T_\lambda f}{\widehat{f}(\lambda)}, g \rangle$ for each $f \in I(M)$, and $\lambda \in \mathcal{S}$, $\Re(\lambda) \neq 0, \frac{1}{2}, 1$, where $T_\lambda f$ is an element of $L^1(G//K)$ such that

$$\widehat{T_\lambda f}(s) = \frac{\widehat{f}(\lambda) - \widehat{f}(s)}{\lambda(1-\lambda) - s(1-s)}.$$

Find this element explicitly.

(5) Use the expressions for $\mathfrak{C}[g]$ to find estimates for the growth of $\mathfrak{C}[g]$ near the boundary of $\mathcal{S}$.

(6) Use maximum modulus techniques to show that $\mathfrak{C}[g] = 0$ in case (a) and that $\mathfrak{C}[g]$ vanishes at $\infty$ and has at most simple poles in $\{0, 1\}$ in case (b).

(7) Show that $\{b_\lambda \mid \lambda \in \mathbb{C} \setminus \mathcal{S}\}$ span a dense subspace of $L^1(G//K)$.

(8) Conclude that $g = 0$ in case (a) and that $g$ is constant in case (b). As $g$ is an arbitrary bounded function annihilating $I(M)$, the theorem follows.

In (1) we find that $b_\lambda(x) = 2Q_{\lambda-1}(x)$, where $Q_{\lambda-1}$ is the Legendre function of the second kind, meets the requirements. This is done using formulae and asymptotics of the Legendre functions.

In (3) $\mathfrak{C}[g]$ is holomorphic outside $\mathcal{S}$ by known properties of Legendre functions. The existence of a holomorphic extension is one of the basic properties of the resolvent transform method and is proved using Banach algebras techniques in the following way.

Let $L^1_\delta(G//K)$ be $L^1(G//K)$ with the unit $\delta$, the dirac measure at the identity of $G$, attached, and assume $g \in L^\infty(G//K)$ annihilates $I(M)$ and that $s \notin \mathcal{S}$. The Fourier transform of the element $\delta - [\lambda(1-\lambda) - s(1-s)]b_s + I(M)$ of $L^1_\delta(G//K)/I(M)$ never vanishes in the maximal ideal space of this quotient algebra ($\{\infty\}$ in case (a) and $\{\infty, 0, 1\}$ in case (b)). This means that it is invertible there and the element $A_\lambda = (\delta - (\lambda(1-\lambda) - s(1-s))b_s + I(M))^{-1}(b_s + I(M))$ of $L^1(G//K)/I(M)$ is well defined. Comparing Fourier transforms, we see that $A_\lambda = b_\lambda + I(M)$ for $\lambda \in \mathbb{C} \setminus \mathcal{S}$. Hence $\mathfrak{C}[g](\lambda) = \langle A_\lambda, g \rangle$ gives the holomorphic extension of $\mathfrak{C}[g]$.

In (4)

$$(T_\lambda f)(t) = P_{\lambda-1}(t) \int_{x>t} f(x) Q_{\lambda-1}(x) \, dx - Q_{\lambda-1}(t) \int_{x>t} f(x) P_{\lambda-1}(x) \, dx,$$

where $t = \cosh(2\zeta)$, meets the requirements. This is done using formulae and asymptotics of the Legendre functions.

In (5) we use Legendre functions formulae to show that $|\mathfrak{C}[g](\lambda)| \leq \frac{C}{\Re(\lambda)\Re(1-\lambda)}$, for $\lambda \notin \mathcal{S}$, and that $|\langle T_\lambda f, g \rangle| \leq \frac{C\|f\|_1}{(\Re(\lambda)\Re(1-\lambda)\Re(\frac{1}{2}-\lambda))^2}$ for $\lambda \in \mathcal{S}$, $\Re(\lambda) \neq 0, \frac{1}{2}, 1$.

In (6) the proof that $\mathfrak{C}[g]$ is zero at infinity is exactly the same as in [10, Theorem 3.3]. To apply this method in the proof that $\mathfrak{C}[g]$ has simple poles in $\{0, 1\}$ in case (b), we use a conformal mapping $\beta$ which "rotates" the strip $\mathcal{S}$ and maps $\infty$ to 0 and conformal mappings which map half strips in the neighborhood of $\infty$ to sectors in the neighborhood of 0. Here is a short description.

The Ahlfors-Heins theorem and the decay conditions in Theorem 1 show that in case (a), $\mathfrak{C}[g](x+iy) = O(\exp(\varepsilon e^{\pi y}))$ as $y \to +\infty$ for every $\varepsilon > 0$ and for almost all $0 < x < 1$ and in case (b), $|z^2\mathfrak{C}[g](z)| = O(\exp(\varepsilon|z|^{-1}))$ as $z \to 0$ for every $\varepsilon > 0$ and for almost all $0 < q < 1$, where $z = \beta(q+it)$.

These estimates are even better than is necessary to apply the Phragmén-Lindelöf theorem for the strip $\mathcal{S}$ (or for the sector $|\arg z| \leq \pi/2$ in case (b)). But there are still three obstacles to overcome before we can use the Phragmén-Lindelöf theorem. First, "almost all" is not enough. Second, the estimates are not uniform in $x$ in case (a) or in $q$ in case (b). Third, $\mathfrak{C}[g]$ is bounded on the boundary of any strip wider than $\mathcal{S}$; however, it may be unbounded on the boundary of $\mathcal{S}$. (The formulation of this problem in case (b) is: $z^2\mathfrak{C}[g](z)$ is bounded on the boundary of any sector wider than the right half plane; however, it may not be bounded on the $y$-axis.) Remember, however, that $|\mathfrak{C}[g](z)| \leq C_1/|\Re(z)\Re(1-z)|$ for $z \notin \mathcal{S}$ in case (a), and that there is a similiar estimate in case (b).

To overcome the first and second problems, we use the Phragmén-Lindelöf theorem on very narrow strips or sectors. For this we need some estimates which are provided by the factorization of $H^\infty$ functions and, close to the boundary, by Levinson's loglog theorem. The implementation of Levinson's theorem in case (b) is different from the one in [10]. Here we need to use Domar's quantitative proof of the theorem as presented in [13, p. 376]. The third problem is solved by use of the Phragmén-Lindelöf argument for domains similiar to strips (or similiar to sectors in case (b)), i.e., by the first distortion inequality of Ahlfors.

Solving the above problems, we get that $\mathfrak{C}[g]$ is bounded in the neighborhood of $\infty$ and that $z^2\mathfrak{C}[g](z)$ is bounded in the neighborhood of $0$. In case (a) we get by Liouville's theorem and the above estimates that $\mathfrak{C}[g] = 0$. In case (b) we get that $\mathfrak{C}[g]$ is zero at $\infty$; and since $z\mathfrak{C}[g]$ is trivially bounded in a sector, it follows that $\mathfrak{C}[g]$ has at most simple poles in $\{0, 1\}$.

In (7), functions whose Fourier transforms decay very fast and extend analytically to a wider strip are dense in $L^1(G//K)$ (see [3, Lemma 1.3]). Each such function can be approximated by linear combination of the $b_\lambda$'s by replacing the integral in the Cauchy representation theorem by Riemann sums.

In (8) the density of the $b_\lambda$'s and (6) give that $g = 0$ in case (a). Let $m \in L^\infty(G//K)$ be the constant function whose value is $\mathrm{Res}(\mathfrak{C}[g](\lambda); 0)$. By (6) $\mathfrak{C}[g] - \mathfrak{C}[m]$ is bounded and vanishes at infinity; hence it is identically zero by Liouville's theorem. Step (7) again implies $g = m$. □

*Proof of Theorem* 2. This follows from Theorem 1 in the same way that the proof of Theorem 3.1 in [3] uses the weaker version of Wiener's theorem proved there. □

*Proof of Corollary* 3. There is an interval $[a, b] \subset (0, \infty)$ where $0 < \mu[a, b] < 1$. Since $P_{x-1}(\cosh(2\zeta)) \leq 1$,

$$\int_0^\infty (1 - P_{x-1}(\cosh(2\zeta)))\, d\mu(\zeta) \geq \int_a^b (1 - P_{x-1}(\cosh(2\zeta)))\, d\mu(\zeta)$$
$$= \mu([a, b]) \int_a^b (1 - P_{x-1}(\cosh(2\zeta))) \frac{d\mu(\zeta)}{\mu([a, b])}.$$

Hence

$$0 \leq -x\log(1-\hat{\mu}(x)) = -x\log(\int_0^\infty (1-P_{x-1}(\cosh(2\zeta)))\,d\mu(\zeta))$$
$$\leq -x\log\mu([a,\,b]) - x\log\int_a^b (1-P_{x-1}(\cosh(2\zeta)))\,\frac{d\mu(\zeta)}{\mu([a,\,b])},$$

and this expression converges to zero using Jensen's inequality, the Lebesgue dominated convergence theorem, and properties of the Legendre functions. The other conditions in Theorem 2 are trivially satisfied. □

*Proof of Theorem* 4. This follows from Theorem 1 in the same way that the proof of Theorem 4.3 in [3] uses the weaker version of Wiener's theorem proved there. □

Y. Ben Natan and Y. Benyamini, Department of mathematics, Technion, Israel Institute of Technology, Haifa, 32000 Israel

*E-mail address*, Y. Ben Natan: `bennatan@techunix.technion.ac.il`

*E-mail address,* Y. Benyamini: `mar29aa@technion`

(H. Hedenmalm) Uppsala Universitet, Matematiska institutionen, Department of Mathematics, Box 480, S-751 06 Uppsala, Sweden
*E-mail address*: `haakan@tdb.uu.se`

(Y. Weit) Department of Mathematics, University of Haifa, Haifa, 31999 Israel
*E-mail address*: `rsma604@haifauvm`